# QRnet: optimal regulator design with LQR-augmented neural networks

Tenavi Nakamura-Zimmerer[1], Qi Gong[1], and Wei Kang[2]

*Abstract*— In this paper we propose a new computational method for designing optimal regulators for high-dimensional nonlinear systems. The proposed approach leverages physics-informed machine learning to solve high-dimensional Hamilton-Jacobi-Bellman equations arising in optimal feedback control. Concretely, we augment linear quadratic regulators with neural networks to handle nonlinearities. We train the augmented models on data generated without discretizing the state space, enabling application to high-dimensional problems. We use the proposed method to design a candidate optimal regulator for an unstable Burgers' equation, and through this example, demonstrate improved robustness and accuracy compared to existing neural network formulations.

## I. INTRODUCTION

WHILE the linear quadratic regulator (LQR) is firmly established as one of the most powerful tools in linear control, the design of optimal regulators for nonlinear systems continues to challenge the control community. The bottleneck in optimal feedback design is the need to solve a Hamilton-Jacobi-Bellman (HJB) partial differential equation (PDE). Due to the well-known "curse-of-dimensionality," this can be extremely difficult for high-dimensional nonlinear systems.

For this reason, there is an extensive literature on methods for approximating solutions of HJB equations. Some key examples include series expansions [1]–[5], level set methods [6], patchy dynamic programming [7], [8], semi-Lagrangian methods [9], [10], method of characteristics and Hopf formula-based algorithms [11], [12], tensor-based methods [13], and polynomial approximation [14]. Unfortunately, many of these methods are limited to moderate dimensions, local solutions, or dynamics with certain algebraic structure.

In recent years, neural networks (NNs) have gained considerable attention as a promising tool for high-dimensional problems since they can avoid the use of spatial grids. Many NN-based methods represent the solution of the HJB equation – called the *value function* – with a NN and minimize the residual of the HJB PDE and boundary conditions at randomly sampled collocation points [15]–[17]. That is, they solve the HJB PDE in the least-squares sense. [18] propose a method for learning a suboptimal policy and approximate value function locally around some nominal trajectories.

Finally, a number of recent works have demonstrated the potential of data-driven methods for HJB. The core idea is to generate data by solving a number of two-point boundary value problems (BVPs) which describe the characteristics of the value function. These BVPs can be solved independently, without a spatial mesh, and in parallel, thus making the algorithm *causality-free*. This property allows the method to be applied to high-dimensional problems. Given BVP data, one then constructs a model of the value function based on this data. In [19], [20] the value function is calculated with sparse grid interpolation, in [21]–[23] supervised learning is used to train a NN model, and in [24] the value function is approximated by sparse polynomial regression. Lastly, [25] consider a related approach which connects forward-backward stochastic differential equations with the HJB equation in stochastic optimal control.

In this paper we propose a physics-informed machine learning method to solve high-dimensional infinite horizon quadratic regularization problems. This important class of control problems arise in the design controllers for regularization or set-point tracking, and has numerous engineering applications in e.g. aerospace, robotics, chemical process control, and distributed parameter systems. The contributions of the present work are to 1) extend the framework introduced in [21]–[23] to infinite horizon problems; and 2) introduce a NN value function model which includes a quadratic term based on the LQR approximation for the linearized dynamics. The proposed model structure mirrors the form of a series expansion with the NN accounting for higher order terms.

We contend that the proposed model structure, which we call *QRnet*, has the following advantages:

- In common practice, linear and nonlinear parts of the control are often treated separately. On the other hand, the QRnet feedback controller smoothly integrates these components to achieve good performance on large domains while retaining the local robustness of LQR.
- Model training is LQR-initialized: rather than learning from scratch, the NN builds on the scaffolding of the LQR value function. As we show in Section IV-B, this can reduce sensitivity to variations in the data set and weight initialization.

[1]Tenavi Nakamura-Zimmerer and Qi Gong are with the Department of Applied Mathematics, Baskin School of Engineering, University of California, Santa Cruz, CA. tenakamu@ucsc.edu, qigong@soe.ucsc.edu
[2]Wei Kang is with the Department of Applied Mathematics, Naval Postgraduate School, Monterey, CA. wkang@nps.edu

## II. PROBLEM SETTING

We consider infinite-horizon nonlinear optimal control problems (OCPs) of the form

$$\begin{cases} \underset{\boldsymbol{u}(\cdot)}{\text{minimize}} & J[\boldsymbol{u}(\cdot)] = \int_0^\infty \mathcal{L}(\boldsymbol{x}, \boldsymbol{u}) dt, \\ \text{subject to} & \dot{\boldsymbol{x}}(t) = \boldsymbol{f}(\boldsymbol{x}, \boldsymbol{u}), \\ & \boldsymbol{x}(0) = \boldsymbol{x}_0. \end{cases} \quad (1)$$

Here $\boldsymbol{x} : [0, \infty) \to \mathbb{R}^n$ is the state, $\boldsymbol{u} : [0, \infty) \to \mathbb{R}^m$ is the control, and $\boldsymbol{f}(\boldsymbol{x}, \boldsymbol{u}) : \mathbb{R}^n \times \mathbb{R}^m \to \mathbb{R}^n$ is a Lipschitz continuous vector field. We consider problems with quadratic running cost,

$$\mathcal{L}(\boldsymbol{x}, \boldsymbol{u}) = (\boldsymbol{x} - \bar{\boldsymbol{x}})^T \boldsymbol{Q} (\boldsymbol{x} - \bar{\boldsymbol{x}}) + (\boldsymbol{u} - \bar{\boldsymbol{u}})^T \boldsymbol{R} (\boldsymbol{u} - \bar{\boldsymbol{u}}), \quad (2)$$

where $\boldsymbol{Q} \in \mathbb{R}^{n \times n}$ is positive semi-definite, $\boldsymbol{R} \in \mathbb{R}^{m \times m}$ is positive definite, and $\bar{\boldsymbol{x}} \in \mathbb{R}^n$, $\bar{\boldsymbol{u}} \in \mathbb{R}^m$ are a (possibly unstable) fixed point of the dynamics such that $\boldsymbol{f}(\bar{\boldsymbol{x}}, \bar{\boldsymbol{u}}) = \boldsymbol{0}$. This standard cost function is a natural choice for regularization or set-point tracking problems where we want to stabilize the objective state $\bar{\boldsymbol{x}}$.

Throughout this paper, we assume that the OCP (1) is well-posed, i.e. an optimal feedback control $\boldsymbol{u}^*(t)$ exists. This ensures that the optimal controlled trajectory, $\boldsymbol{x}^*(t)$, satisfies $\lim_{t \to \infty} \mathcal{L}(\boldsymbol{x}^*(t), \boldsymbol{u}^*(t)) = 0$. Due to real-time application requirements, we would like to design a closed-loop feedback controller, $\boldsymbol{u}^*(t) = \boldsymbol{u}^*(\boldsymbol{x}(t))$, which can be evaluated online given any measurement of $\boldsymbol{x}$.

### A. The Hamilton-Jacobi-Bellman equation

Following the standard procedure in optimal control, to compute the optimal feedback control we begin by defining the value function $V : \mathbb{R}^n \to \mathbb{R}$ as the optimal cost-to-go of (1) starting at the point $\boldsymbol{x}(0) = \boldsymbol{x}$. That is,

$$V(\boldsymbol{x}) := J[\boldsymbol{u}^*(\cdot)]. \quad (3)$$

It can be shown that the value function is the unique viscosity solution [26] of the steady state HJB PDE,

$$\begin{cases} \min_{\boldsymbol{u}} \{ \mathcal{L}(\boldsymbol{x}, \boldsymbol{u}) + V_{\boldsymbol{x}}^T(\boldsymbol{x}) \boldsymbol{f}(\boldsymbol{x}, \boldsymbol{u}) \} = 0, \\ V(\bar{\boldsymbol{x}}) = 0, \end{cases} \quad (4)$$

where we denote $V_{\boldsymbol{x}} := \partial V / \partial \boldsymbol{x}$. If (4) can be solved (in the viscosity sense), then it provides both necessary and sufficient conditions for optimality.

Given the value function $V(\cdot)$, we define the Hamiltonian

$$\mathcal{H}(\boldsymbol{x}, V_{\boldsymbol{x}}, \boldsymbol{u}) := \mathcal{L}(\boldsymbol{x}, \boldsymbol{u}) + V_{\boldsymbol{x}}^T \boldsymbol{f}(\boldsymbol{x}, \boldsymbol{u}). \quad (5)$$

The optimal control satisfies the Hamiltonian minimization condition,

$$\boldsymbol{u}^*(\boldsymbol{x}) = \boldsymbol{u}^*(\boldsymbol{x}; V_{\boldsymbol{x}}(\boldsymbol{x})) = \arg\min_{\boldsymbol{u}} \mathcal{H}(\boldsymbol{x}, V_{\boldsymbol{x}}, \boldsymbol{u}). \quad (6)$$

Then if we can solve (4), the optimal feedback control is obtained online as the solution of (6).

But as discussed in Section I, solving (4) is extremely challenging. Thus a common approach is to linearize the dynamics about $(\boldsymbol{x} = \bar{\boldsymbol{x}}, \boldsymbol{u} = \bar{\boldsymbol{u}})$ to obtain the linear system

$$\begin{cases} \frac{d}{dt}(\boldsymbol{x} - \bar{\boldsymbol{x}}) \approx \boldsymbol{A}(\boldsymbol{x} - \bar{\boldsymbol{x}}) + \boldsymbol{B}(\boldsymbol{u} - \bar{\boldsymbol{u}}), \\ \boldsymbol{A} := \frac{\partial \boldsymbol{f}}{\partial \boldsymbol{x}}(\bar{\boldsymbol{x}}, \bar{\boldsymbol{u}}), \quad \boldsymbol{B} := \frac{\partial \boldsymbol{f}}{\partial \boldsymbol{u}}(\bar{\boldsymbol{x}}, \bar{\boldsymbol{u}}). \end{cases} \quad (7)$$

Under the mild conditions that $(\boldsymbol{A}, \boldsymbol{B})$ is controllable and $(\boldsymbol{A}, \boldsymbol{Q}^{1/2})$ is observable, the value function of the OCP with the linear dynamics (7) and quadratic cost (2) is

$$V^{\text{LQR}}(\boldsymbol{x}) = (\boldsymbol{x} - \bar{\boldsymbol{x}})^T \boldsymbol{P} (\boldsymbol{x} - \bar{\boldsymbol{x}}), \quad (8)$$

where $\boldsymbol{P} \in \mathbb{R}^{n \times n}$ is a positive semi-definite matrix satisfying the Riccati equation,

$$\boldsymbol{Q} + \boldsymbol{A}^T \boldsymbol{P} + \boldsymbol{P} \boldsymbol{A} - \boldsymbol{P} \boldsymbol{B} \boldsymbol{R}^{-1} \boldsymbol{B}^T \boldsymbol{P} = \boldsymbol{0}. \quad (9)$$

Furthermore, the resulting state feedback controller is linear with constant gain:

$$\boldsymbol{u}^{\text{LQR}}(\boldsymbol{x}) = -\boldsymbol{K}\boldsymbol{x}, \qquad \boldsymbol{K} = \boldsymbol{R}^{-1} \boldsymbol{B}^T \boldsymbol{P}. \quad (10)$$

This approach has yielded many successful engineering applications, but it is suboptimal and in some cases even fails to stabilize the nonlinear dynamics. For this reason we are interested in finding the value function $V(\cdot)$ for the full nonlinear dynamics.

### B. Pontryagin's Minimum Principle

To make use of (6) for general nonlinear systems, we need an efficient way to approximate the value function and its gradient. Like [19]–[24], rather than solve the full HJB equation (4) directly, we exploit the fact its characteristics evolve according to a two-point BVP, well-known in optimal control as Pontryagin's Minimum Principle (PMP, [27]):

$$\lim_{t_f \to \infty} \begin{cases} \dot{\boldsymbol{x}}(t) = \mathcal{H}_{\boldsymbol{\lambda}} = \boldsymbol{f}(\boldsymbol{x}, \boldsymbol{u}^*(\boldsymbol{x}; \boldsymbol{\lambda})), & \boldsymbol{x}(0) = \boldsymbol{x}_0, \\ \dot{\boldsymbol{\lambda}}(t) = -\mathcal{H}_{\boldsymbol{x}}(\boldsymbol{x}, \boldsymbol{\lambda}, \boldsymbol{u}^*(\boldsymbol{x}; \boldsymbol{\lambda})), & \boldsymbol{\lambda}(t_f) = \boldsymbol{0}, \\ \dot{v}(t) = -\mathcal{L}(\boldsymbol{x}, \boldsymbol{u}^*(\boldsymbol{x}; \boldsymbol{\lambda})), & v(t_f) = 0. \end{cases} \quad (11)$$

Here $\boldsymbol{\lambda} : [0, \infty) \to \mathbb{R}^n$ is called the *costate*. The two-point BVP (11) provides a necessary condition for optimality, and if we further assume that the solution is optimal, then along the characteristic $\boldsymbol{x} = \boldsymbol{x}^*(t; \boldsymbol{x}_0)$ we have

$$V(\boldsymbol{x}) = v(t), \quad V_{\boldsymbol{x}}(\boldsymbol{x}) = \boldsymbol{\lambda}(t), \quad \boldsymbol{u}^*(\boldsymbol{x}) = \boldsymbol{u}^*(t). \quad (12)$$

In general, the BVP (11) admits multiple solutions. So while the characteristics of the value function satisfy (11), there may be other solutions to these equations which are suboptimal and thus not characteristics. In certain problems the characteristics can also intersect, giving rise to non-smooth value functions and difficulties in applying (12).

Optimality of solutions to (11) can be guaranteed under some convexity conditions (see e.g. [28]). For most dynamical systems it is difficult to verify such conditions globally, but we can guarantee optimality locally around an equilibrium point [2]. Addressing the challenge of global optimality is beyond the scope of the present work, so in this paper we assume that solutions of (11) are optimal. Under this assumption, the relationship between PMP and the value function given in (12) holds everywhere.

Note the proposed method can still be applied even when this assumption cannot be verified. In such cases PMP is the prevailing tool for finding *candidate optimal* solutions, and from these the proposed method yields a stabilizing feedback controller which satisfies necessary conditions for optimality.

## III. Neural network value function modeling

### A. Data generation

Like [22]–[24], we generate data by solving the two-point BVP (11) for a set of randomly sampled initial conditions. Critically, these BVPs can be solved independently without knowledge of nearby solutions. Methods based on this idea are referred to as *causality-free* [19], [20], [22]–[24]. Note that the related method proposed in [21] differs slightly because it does not allow one to choose initial conditions freely. For a survey of data generation methods we refer the reader to [29].

*1) Generating infinite-horizon data:* Different from these prior works which consider fixed finite time OCPs, in this paper we are interested in infinite-horizon problems. Notice that the infinite-horizon PMP (11) is obtained with the limit $t_f \to \infty$ of a finite-horizon problem [27]. To reflect this, we solve (11) up to some fixed final time $t_f$ and check if the running cost $\mathcal{L}(x(t_f), u(t_f))$ is smaller than a desired tolerance. If not, we extend the time horizon and – using the previous solution as an initial guess – re-solve the BVP until the running cost is sufficiently small.

Once the running cost is small enough, it follows that the finite-horizon solution approximates the solution of the infinite-horizon problem as further integration should not change the cost significantly. This conclusion is reasonable as closed-loop stability is a necessary condition for finiteness of the value function. Then by applying (12) at each point along the trajectory $x = x(t; x_0)$ and aggregating data from all infinite-horizon BVP solutions, we obtain a data set

$$\mathcal{D} = \left\{ x^{(i)}, V^{(i)}, \lambda^{(i)} \right\}_{i=1}^{|\mathcal{D}|},$$

where $V^{(i)} := V\left(x^{(i)}\right)$, $\lambda^{(i)} := V_x\left(x^{(i)}\right)$, and $|\mathcal{D}|$ denotes the number of data points in $\mathcal{D}$. Note that there is no need to distinguish data from different trajectories as the value function and its gradient are time-independent.

While generating data in this way is efficient because we extract a lot of data from each successful BVP solution, it has the side effect of concentrating a large amount of data near the equilibrium. On the other hand, we are interested in designing controllers which are effective over large regions of the state space and consequently we need data sets which support learning far from the equilibrium. To this end, instead of including the whole trajectory in the data set, we only take points with $t \leq t_f/T$ for a parameter $T \geq 1$. $T$ is chosen to balance efficiency in data generation with the competing objective of adequately representing the entire state space; in this paper we set $T = 3$.

*2) LQR warm start for reliable BVP solution:* In this paper, we solve the two-point BVP using the SciPy [30] implementation of the BVP solver introduced in [31]. This algorithm is highly accurate but convergence is highly sensitive to the initial guess for $x(t)$ and $\lambda(t)$. Furthermore, convergence is increasingly dependent on good initializations as we increase $t_f$ to approximate the infinite-horizon problem.

To mitigate this difficulty we simulate the dynamics up to some large final time $t_f$ with an LQR controller (10) to close the loop. This provides a guess for the optimal state trajectory, and a guess for the costate can be obtained with the LQR approximation $\lambda(t) \approx 2Px(t)$. While the costate guess is often far from perfect, we find that it is usually close enough to facilitate reliable convergence over large time horizons. We refer to this strategy as *LQR warm-start*.

### B. Neural network architecture

We employ a simple and intuitive architecture to model the value function. The main idea is to augment an LQR value function approximation with a NN which accounts for nonlinearities. The LQR value function is computed with respect to the dynamics linearized around $(\bar{x}, \bar{u})$, and provides a good local approximation. The NN corrects and extends the approximation throughout the training domain.

As for the NN, we use a standard fully-connected feed-forward architecture. We denote the output of the network as $W^{\text{NN}}(\cdot)$. Feedforward NNs approximate complicated nonlinear functions by a composition of simpler functions, namely

$$W^{\text{NN}}(x) = g_L \circ g_{L-1} \circ \cdots \circ g_1(x), \quad (13)$$

where the $\ell$th layer, $\ell = 1, \ldots, L$, is defined as $g_\ell(z) = \sigma_\ell(W_\ell z + b_\ell)$. Within each layer, $W_\ell$ and $b_\ell$ are the weight matrices and bias vectors, respectively, and $\sigma_\ell(\cdot)$ denotes a nonlinear *activation function* applied component-wise to its argument. We keep the scalar output layer $g_L(\cdot)$ linear, so $\sigma_L(\cdot)$ is the identity function.

We combine the raw NN prediction (13) with the LQR value function (8) for the linearized dynamics (7) as

$$V^{\text{NN}}(x) = \frac{1}{c} \log\left[1 + cV^{\text{LQR}}(x)\right] + W^{\text{NN}}(x), \quad (14)$$

with a trainable parameter $c > 0$. Intuitively, LQR provides a good approximation near $\bar{x}$. There $V^{\text{LQR}}(x)$ is small and hence $c^{-1} \log\left[1 + cV^{\text{LQR}}(x)\right] \approx V^{\text{LQR}}(x)$ for all $c \in (0, \infty)$. Further away from $\bar{x}$, we have $c^{-1} \log\left[1 + cV^{\text{LQR}}(x)\right] \ll V^{\text{LQR}}(x)$, thereby increasing the relative importance of the corrective NN. The parameter $c$ governs the radius in which this term approximates $V^{\text{LQR}}(x)$; in particular $\lim_{c \to 0} c^{-1} \log\left[1 + cV^{\text{LQR}}(x)\right] = V^{\text{LQR}}(x)$. Notice that the model structure (14) is similar to a series expansion, except that we explicitly reduce the impact of lower order terms away from the linearization point.

Finally, the NN-based feedback control is evaluated by substituting $V_x^{\text{NN}}(\cdot)$ into (6) in place of the gradient of the true value function:

$$u^{\text{NN}}(x) := u^*\left(x; V_x^{\text{NN}}(x)\right). \quad (15)$$

It should be emphasized that the gradient $V_x^{\text{NN}}(\cdot)$ is calculated using automatic differentiation, and is therefore *exact* and computationally efficient.

### C. Physics-informed learning

Suppose we have generated a data set $\mathcal{D}$ as discussed in Section III-A. This data takes the form of input-output pairs: $x^{(i)}$ are the inputs and $\left(V^{(i)}, \lambda^{(i)}\right)$ are the outputs to be modeled. Let $\theta$ denote the collection of model parameters:

$$\theta := \{c\} \cup \{W_\ell, b_\ell\}_{\ell=1}^L.$$

In [21]–[23], the NN is trained by solving a supervised learning problem in which one minimizes the mean square regression loss,

$$\operatorname*{loss}_V(\boldsymbol{\theta}) := \frac{1}{|\mathcal{D}|} \sum_{i=1}^{|\mathcal{D}|} \left[ V^{(i)} - V^{\mathrm{NN}}\left(\boldsymbol{x}^{(i)}; \boldsymbol{\theta}\right) \right]^2, \quad (16)$$

plus a gradient regularization term

$$\operatorname*{loss}_{\boldsymbol{\lambda}}(\boldsymbol{\theta}) := \frac{1}{|\mathcal{D}|} \sum_{i=1}^{|\mathcal{D}|} \left\| \boldsymbol{\lambda}^{(i)} - V_{\boldsymbol{x}}^{\mathrm{NN}}\left(\boldsymbol{x}^{(i)}; \boldsymbol{\theta}\right) \right\|^2. \quad (17)$$

This serves as a form of *physics-informed* regularization. The term "physics-informed" is borrowed from [32] which partially inspired the authors' previous work [22], [23]. By incorporating the prior knowledge that $\boldsymbol{\lambda}(t) = V_{\boldsymbol{x}}(\boldsymbol{x}(t))$, we maximize the information extracted from the available data and obtain more optimal feedback laws [21]–[23].

In this work we impose an additional penalty on deviating from the optimal control:

$$\operatorname*{loss}_{\boldsymbol{u}}(\boldsymbol{\theta}) := \frac{1}{|\mathcal{D}|} \sum_{i=1}^{|\mathcal{D}|} \left\| \boldsymbol{u}^*\left(\boldsymbol{x}^{(i)}\right) - \boldsymbol{u}^{\mathrm{NN}}\left(\boldsymbol{x}^{(i)}; \boldsymbol{\theta}\right) \right\|^2. \quad (18)$$

Minimizing this control penalty term contributes directly toward the ultimate goal of using the NN for optimal feedback by effectively enforcing the Hamiltonian minimization condition (6) on the learned feedback policy. We now arrive at the following optimization problem which we solve to train the NN:

$$\operatorname*{minimize}_{\boldsymbol{\theta}} \quad \operatorname*{loss}(\boldsymbol{\theta}) := \operatorname*{loss}_V(\boldsymbol{\theta}) + \mu_{\boldsymbol{\lambda}} \operatorname*{loss}_{\boldsymbol{\lambda}}(\boldsymbol{\theta}) + \mu_{\boldsymbol{u}} \operatorname*{loss}_{\boldsymbol{u}}(\boldsymbol{\theta}), \quad (19)$$

where $\mu_{\boldsymbol{\lambda}}, \mu_{\boldsymbol{u}} \geq 0$ are scalar weights.

To quantify the accuracy of the model, we generate two data sets from *independently drawn* initial conditions. During training, the network observes only data points from the training set $\mathcal{D}_{\mathrm{train}}$. The other data set, which we call the validation set $\mathcal{D}_{\mathrm{val}}$, is reserved for evaluating the NN accuracy *after* training. Good validation performance indicates that the NN generalizes well, i.e. it did not overfit the training data. The ability to empirically measure model accuracy in this way is a key feature of causality-free methods.

## IV. APPLICATION TO DISTRIBUTED PARAMETER SYSTEM

In this section we explore the effectiveness of the proposed algorithm by solving a 64-dimensional OCP arising from a Chebyshev pseudospectral (PS) discretization of a modified Burgers' equation with a destabilizing reaction term. Stabilization of Burgers' equation is a common benchmark problem in distributed parameter systems and similar problems have recently been considered in e.g. [5], [14], [22], [23].

Let $X(t,\xi) : [0,\infty) \times [-1,1] \to \mathbb{R}$ satisfy the following one-dimensional controlled PDE with Dirichlet boundary conditions:

$$\begin{cases} X_t = -\frac{1}{2}(X^2)_\xi + \nu X_{\xi\xi} + \alpha(\xi) X e^{-\beta X} + \boldsymbol{b}^T(\xi) \boldsymbol{u}(t) \\ \qquad\qquad\qquad \text{for } t > 0, \xi \in (-1,1), \\ X(t,-1) = X(t,1) \quad \text{for } t > 0, \\ X(0,\xi) = X_0(\xi) \quad \text{for } \xi \in (-1,1). \end{cases} \quad (20)$$

Here $\nu, \beta > 0$ are scalar parameters, $\alpha : (-1,1) \to \mathbb{R}$, and $\boldsymbol{b} : (-1,1) \to \mathbb{R}^m$. The control $\boldsymbol{u} : [0,\infty) \to \mathbb{R}^m$ is designed to stabilize the *open-loop unstable* origin by solving the PDE-constrained OCP

$$\begin{cases} \operatorname*{min.}_{\boldsymbol{u}(\cdot)} \quad J[\boldsymbol{u}(\cdot)] = \int_0^\infty \left( \|X\|^2_{L^2_{(-1,1)}} + R\, \boldsymbol{u}^T \boldsymbol{u} \right) dt, \\ \text{s.t.} \quad \text{Eq. (20)}, \end{cases} \quad (21)$$

where

$$\|X\|^2_{L^2_{(-1,1)}} := \int_{-1}^1 X^2(t,\xi)\, d\xi.$$

We define

$$\alpha(\xi) = \begin{cases} -\kappa \left(\xi + \frac{1}{5}\right)\left(\xi - \frac{1}{5}\right), & \xi \in \left[-\frac{1}{5}, \frac{1}{5}\right], \\ 0, & \xi \notin \left[-\frac{1}{5}, \frac{1}{5}\right]. \end{cases}$$

and consider the case with $m = 2$ actuators active on compact supports defined by

$$\boldsymbol{b}(\xi) = \begin{pmatrix} \begin{cases} -\kappa\left(\xi + \frac{4}{5}\right)\left(\xi + \frac{2}{5}\right), & \xi \in \left[-\frac{4}{5}, -\frac{2}{5}\right], \\ 0, & \xi \notin \left[-\frac{4}{5}, -\frac{2}{5}\right], \end{cases} \\ \begin{cases} -\kappa\left(\xi - \frac{2}{5}\right)\left(\xi - \frac{4}{5}\right), & \xi \in \left[\frac{2}{5}, \frac{4}{5}\right], \\ 0, & \xi \notin \left[\frac{2}{5}, \frac{4}{5}\right]. \end{cases} \end{pmatrix}$$

We set $\nu = 0.02$, $\beta = 0.1$, $\kappa = 25$, and $R = 0.5$.

### A. Pseudospectral discretization

To solve (21) using our framework, we perform Chebyshev PS collocation to transform (20) into a system of $n$ ordinary differential equations (ODEs). Following [33], let

$$\xi_j = \cos(j\pi/n), \qquad j = 0, 1, \ldots n, n+1.$$

We collocate $X(t,\xi)$ at the non-boundary nodes $\xi_j$, $j = 1, 2, \ldots, n$, and set $X(t, \xi_0) = X(t, \xi_{n+1}) = 0$ to account for the boundary conditions. We then define

$$\boldsymbol{x}(t) := \begin{pmatrix} X(t, \xi_1), & X(t, \xi_2), & \ldots, & X(t, \xi_n) \end{pmatrix}^T$$

and construct Chebyshev differentiation matrices $\boldsymbol{D}, \boldsymbol{D}^2 \in \mathbb{R}^{n \times n}$. Hence the PDE (20) becomes

$$\dot{\boldsymbol{x}} = -\frac{1}{2} \boldsymbol{D} \boldsymbol{x}^2 + \nu \boldsymbol{D}^2 \boldsymbol{x} + \boldsymbol{\alpha} \odot \boldsymbol{x} \odot e^{-\beta \boldsymbol{x}} + \boldsymbol{B} \boldsymbol{u}, \quad (22)$$

where $\boldsymbol{x}^2 := \boldsymbol{x} \odot \boldsymbol{x}$, "$\odot$" denotes elementwise multiplication, and $\boldsymbol{\alpha} \in \mathbb{R}^n$ and $\boldsymbol{B} \in \mathbb{R}^{n \times m}$ are the collocated versions of $\alpha(\xi)$ and $\boldsymbol{b}^T(\xi)$.

The integral appearing in the cost function is conveniently approximated by Clenshaw-Curtis quadrature [33]. Let $w_j$, $j = 1, 2, \ldots, n$ be the non-boundary Clenshaw-Curtis quadrature weights. Then

$$\begin{cases} \|X\|^2_{L^2_{(-1,1)}} = \int_{-1}^1 X^2(t,\xi)\, d\xi \approx \boldsymbol{x}^T \boldsymbol{Q} \boldsymbol{x}, \\ \boldsymbol{Q} = \operatorname{diag}(w_1, \quad w_2, \quad \ldots, \quad w_n). \end{cases}$$

Now the original OCP (21) can be reformulated as a quadratic cost ODE-constrained problem,

$$\begin{cases} \operatorname*{min.}_{\boldsymbol{u}(\cdot)} \quad J_n[\boldsymbol{u}(\cdot)] = \int_0^\infty \left( \boldsymbol{x}^T \boldsymbol{Q} \boldsymbol{x} + R \boldsymbol{u}^T \boldsymbol{u} \right) dt, \\ \text{s.t.} \quad \text{Eq. (22)}. \end{cases} \quad (23)$$

## B. Numerical results

We now present results of our method for solving the OCP (23) with $n = 64$ collocation points. The algorithm can handle higher dimensions, but we find that $n = 64$ points are enough to resolve the stiff dynamics for all initial conditions tested. We compare the accuracy of QRnet with that of a standard NN trained using supervised learning (comparable to [21]–[23]), as well as with a straight LQR approximation.

We consider initial conditions which are sums of sine functions with uniform random coefficients:

$$X_0(\xi) = \sum_{k=1}^{10} a_k \sin(k\pi\xi), \quad a_k \sim \mathcal{U}(-1/k, 1/k), \quad (24)$$

and use the following hyperparameters for training:
- in the NN component $W^{\text{NN}}(\cdot)$ we use $L = 5$ layers each with 32 neurons, and apply the $\tanh(\cdot)$ activation function to all hidden layers $\ell = 1, \ldots, 4$;
- we set the weights in the loss function (19) to $\mu_{\boldsymbol{\lambda}} = 0$ and $\mu_{\boldsymbol{u}} = 5$, and optimize using L-BFGS-B [34];
- we implement the model in TensorFlow 1.11 [35] and train it on an NVIDIA RTX 2080Ti GPU.

First we study the sensitivity of the method with respect to variations in the data set and parameter initialization. This is important as collecting data can be time-consuming and NN training is a highly non-convex optimization problem. To this end we vary the number of trajectories in the training set $\mathcal{D}_{\text{train}}$, and for each different data set size, conduct ten trials with different randomly generated training trajectories and NN weight initializations. For validation we build a data set of 400 trajectories totaling $|\mathcal{D}_{\text{val}}| = 66060$ data points.

From the results shown in Fig. 1, it is clear that both the plain NN and QRnet vastly outperform LQR for value function reconstruction and optimal control prediction. Furthermore, the distribution of validation errors skews lower for QRnet than the plain NN. This suggests that QRnet is less sensitive to variations in the data set and parameter initialization, and thus enables the use of smaller data sets. The training time for smaller data sets is also shorter. Together, these properties facilitate rapid prototyping over an iterative design process, and combine naturally with the progressive model refinement strategy proposed in [22].

Next we select NN and QRnet models trained on the same set of 64 trajectories with nearly identical validation accuracy. We compare the performance of the LQR, NN, and QRnet feedback controllers to the open-loop optimal control, across 1200 simulations for initial conditions of different size, $\|X_0\|_{L^2_{(-1,1)}} = 0.1, 0.2, \ldots, 1.2$. Results are shown in Fig. 2. As expected, the performance of LQR steadily degrades away from the origin, while the NN and QRnet controllers do well throughout the training domain. The latter two appear to be largely equivalent, except for a few initial conditions where the NN accumulates significant additional cost while LQR and QRnet do not. This supports the idea that QRnet inherits some local robustness from LQR, yielding a more reliable controller than a plain NN.

We conclude with a simulation of the dynamics (22) with QRnet feedback, for a typical random initial condition. Results are plotted in Fig. 3. There we can see that QRnet stabilizes the system and closely approximates the optimal control, even where LQR deviates from it.

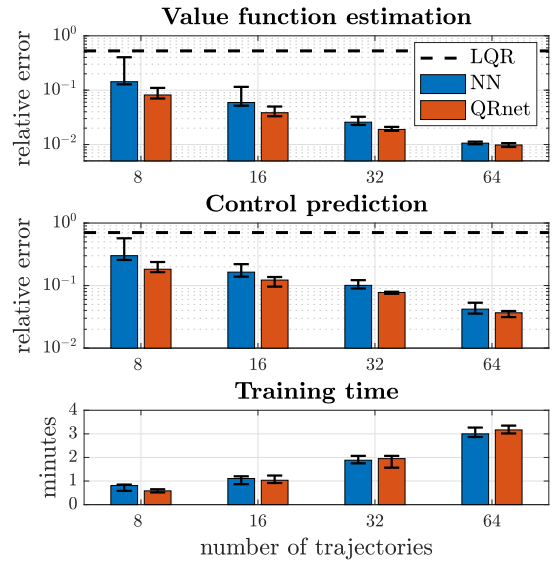

Fig. 1: Relative mean absolute error in estimating the value, relative mean $L^2$ error in predicting the optimal control, and training time, depending on the number of trajectories seen during training. The bar graph height shows the median over ten trials and error bars indicate the 25th and 75th percentiles.

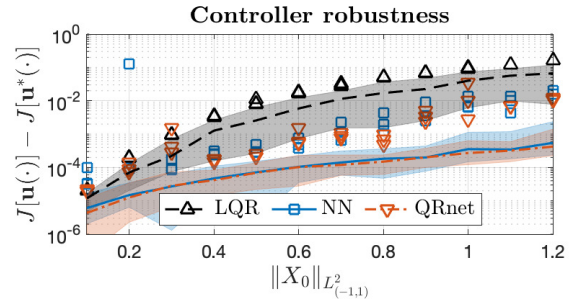

Fig. 2: Difference between controlled and optimal costs, depending on the norm of the initial condition. Shaded areas cover the 15th to 85th percentiles, lines show medians, and symbols pick out the top three outliers for each group of initial conditions.

## V. CONCLUDING REMARKS

In this paper we have presented *QRnet*, an extension of the causality-free physics-informed learning framework [22], [23] to infinite-horizon OCPs. This extension is comprised of two main features: efficient infinite-horizon data generation and a structurally-motivated NN architecture. By way of the Burgers' benchmark problem, we have illustrated the potential for use in high-dimensional nonlinear systems and the improved robustness and performance the LQR-augmented model architecture.

Much remains to be explored in the development of deep learning approaches for feedback design. For instance, we are interested in applying the proposed method to problems with control constraints and non-quadratic costs with locally quadratic approximations. It will also be useful to study how beneficial the model architecture is, depending on how well the original problem can be approximated by a linear quadratic one. Finally, since the computational framework depends crucially on data generation, in future work we plan to study different strategies for improving the robustness and efficiency of this step.

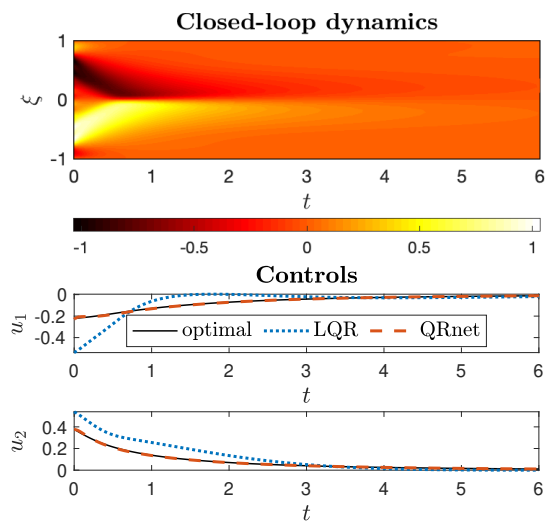

Fig. 3: Closed-loop dynamics and controls of the Burgers' system (22) for a random initial condition. The full simulation is over $t \in [0, 30]$ but we show only $t \in [0, 6]$.